\documentclass{amsart}
\newtheorem{theorem}{Theorem}

\newtheorem{corollary}[theorem]{Corollary}

\newtheorem{definition}[theorem]{Definition}
\newtheorem{example}[theorem]{Example}

\newtheorem{lemma}[theorem]{Lemma}
\newtheorem{proposition}[theorem]{Proposition}

\numberwithin{equation}{section}

\renewcommand{\oddsidemargin}{0.5cm}

\def\C{\mathbb C}
\def\R{\mathbb R}
\def\X{\mathbb X}

\def\Z{\mathbb Z}
\def\Y{\mathbb Y}
\def\Z{\mathbb Z}
\def\N{\mathbb N}
\def\cal{\mathcal}

\begin{document}
\title{Massera Type Theorems for Abstract Functional Differential Equations}
\author{Qing Liu}
\address{School of Mathematical Science, Beijing Normal University, Beijing 100875, People's Republic of China }
\email{bnumlq@mail.bnu.edu.cn}
\author{Nguyen Van Minh}
\address{Department of Mathematics, University of West Georgia, Carrollton, GA 30118}
\email{vnguyen@westga.edu}
\author{G. Nguerekata}
\address{Department of Mathematics, Morgan State University,
1700E. Cold Spring Lane, Baltimore, MD21251}
\email{gnguerek@jewel.morgan.edu}
\author{Rong Yuan}
\address{School of Mathematical Science, Beijing Normal University, Beijing 100875, People's Republic of China }
\email{ryuan@bnu.edu.cn}
\thanks{The third author N.V.M. thanks the Beijing Normal University for the support during his visit in summer 2006. The authors
of the paper are grateful to anonymous referees and the editor for
pointing out several inaccuracies and for the suggestions to improve
the presentation of the paper.}

\date{\today}
\begin{abstract}
The paper is concerned with conditions for the existence of almost
periodic solutions of the following abstract
 functional differential equation
$ \dot u(t) = Au(t) + [{\cal B}u](t) +f(t), $ where $A$ is a closed
operator in a Banach space $\X$, $\cal B$ is a general bounded
linear operator in the function space of all $\X$-valued bounded and
uniformly continuous functions that satisfies a so-called {\it
autonomous} condition. We develop a general procedure to carry out
the decomposition that does not need the well-posedness of the
equations. The obtained conditions are of Massera type, which are
stated in terms of spectral conditions of the operator ${\cal
A}+{\cal B}$ and the spectrum of $f$. Moreover, we give conditions
for the equation not to have quasi-periodic solutions with different
structures of spectrum. The obtained results extend previous ones.
\end{abstract}

\maketitle

\section{Introduction}
In this paper we are concerned with the existence of almost periodic
solutions to abstract functional differential equations of the form
\begin{equation}\label{1}
\dot u(t) = Au(t) + [{\cal B}u](t) +f(t),
\end{equation}
where $A$ is a closed operator on a Banach space $\X$, ${\cal B}$
is a bounded linear {\it autonomous} operator on $BUC(\R,\X)$, $f$
is a $\X$-valued almost periodic function.

As is well known, in \cite{mas} Massera studied the following
linear ordinary differential equation
\begin{equation}\label{2}
\dot u(t) = A(t)u(t) + f(t),
\end{equation}
where $A$ and $f$ are continuous, periodic with the same period
$\tau$, and proved a classical theorem that is often referred to
as Massera Theorem, saying that Eq.(\ref{2}) has a periodic
solution with period $\tau$ if and only if it has a bounded
solution on the positive half line. The Massera Theorem has been
extended to various kinds of differential equations. We refer the
reader to \cite{chohal,lilinli,liconlinliu,shinai,benezz,her} and
the references therein for more information in this direction. The
method employed in these works is to prove the existence of
periodic solutions via the existence of fixed points of the
associated period maps. This method is no longer valid for more
general problems dealing with almost periodic solutions, as one
has no period map associated with a non-periodic equation. A new
idea to study Massera type theorems for almost periodic solutions
was developed in \cite{naiminshi} that is referred to as spectral
decomposition of bounded solutions. This idea can be carried out
via the concept of evolution semigroups. Extensions of the results
in \cite{naiminshi} to various kind of equations were obtained in
\cite{furnaimin,minmin,murnaimin2,ezzjaz}. Note that the methods
used in these papers are more or less based on the well-posedness
of the equations, so they are no longer valid for the general
setting of this paper in which the operator $A$ may not generate a
semigroup, and especially, $\cal B$ is a very general functional
operator.

\medskip In this paper we will develop a transparent operator theoretical
framework for decomposition of bounded solutions of Eq. (\ref{1})
with very general assumptions on $A$ and $\cal B$. This framework
does not require any conditions for the well-posedness, for the
variation-of-constants formula, so it can be applied to many kind
of functional equations. We note that most of periodic functional
evolution equations considered in previous works (see e.g.
\cite{furnaimin,murnaimin2,shinai}) on Massera type theorems can
be reduced to equations with constant coefficients by the partial
Floquet representation developed in \cite{hen}. Therefore, our
results obtained in this paper can be extended to these equations.

\medskip
We now briefly outline our paper. In the next section we review
some concepts such as almost periodicity, spectrum of a bounded
function, and especially, the Loomis Theorem in the infinite
dimensional case that is a key for our results. The section that
follows contains the main results of this paper with two Massera
type theorems (Theorems \ref{th 3}, \ref{th 8}) whose refinements
are Corollary \ref{cor 11}. A result on the non-existence of
quasi-periodic solutions is Theorem \ref{cor 12}. Finally, we give
an application to show that our results can be applied to
functional evolution equations of mixed type which are, in
general, not well-posed, that is, the associated Cauchy Problem
may have no solutions.

\section{Preliminaries}

In this section, we first give some notations and definitions. After
recalling the concept of almost periodic functions we present the
spectral theory of functions and several important properties which
we need in the next sections. Finally, two lemmas on decomposition
are proved.

\subsection{Notations and Definitions}

Throughout this paper we will use the following notations: $\N, \Z,
\mathbb{Q}, \R, \C$ stand for the sets of natural, integer,
rational, real and complex numbers, respectively. For any complex
number $\lambda$, the notation $Re\lambda$ stands for its real part.
We always denote by $\X$ a given complex Banach space, and
$BC(\R,\X), BUC(\R,\X), AP(\X)$ the spaces of all $\X$-valued
bounded continuous, bounded uniformly continuous and almost periodic
functions on $\R$ with sup-norm, respectively. The translation group
on $BC(\R, \X)$ is denoted by $(S(t))_{t\in \R}$ and its generator
$\mathcal{D} = \frac{d}{dt}$. The spaces $BC(\R,\X), BUC(\R,\X),
AP(\X)$ are translation-invariant. Given two complex Banach spaces
$\X, \Y$, for any linear operator $T$ from $\X$ to $\Y$, as usual
$\sigma(T), \rho(T)$, $R(\lambda, T)=(\lambda- T)^{-1}$ are the
notations of the spectrum, resolvent set, resolvent of $T$.

Let us denote by ${\cal A}$ the operator of multiplication by $A$ on
$BUC(\R,\X)$, that is, an operator with domain
$$D({\cal A}):=\{ g\in BUC(\R,\X):\ g(t)\in D(A), \ \forall t\in \R, \ \mbox{and} \ Ag(\cdot )\in BUC(\R,\X) \} $$
and $[{\cal A}u](t)=Au(t),$ for all $t\in\R$ and $u\in D({\cal A})$.

\medskip
Recall that two closed operators $A$ and $B$ on a Banach space
$\Y$ with non-empty resolvent sets commute in the sense of
resolvent commuting if $R(\lambda, A)R(\mu, B)=R(\mu, B)R(\lambda,
A)$ for all $\lambda\in \rho(A), \mu\in \rho(B).$ Some elementary
properties of commuting operators will be used in this paper whose
proofs are left to the reader. For example, if $A$ generates a
semigroup $(T(t))_{t\ge 0}$, then a closed operator $B$ commutes
with $A$ if and only if it commutes with $T(t)$ for all $t\ge 0$.

\begin{definition}
A bounded linear operator ${\cal B}$ in $BUC(\R,\X)$ is said to be
autonomous if it commutes with the translation group $(S(t))_{t\in
\R}$.
\end{definition}


\subsection{Almost Periodic Functions}
Let $g$ be a continuous function on $\R$ taking values in a complex
Banach space $\X$. $g$ is said to be almost periodic in the sense of
Bohr if to every $\epsilon > 0$, the set
\begin{equation*}
T(g, \epsilon):=\{\tau: \sup_{t\in
\R}\|g(t+\tau)-g(t)\|\leq\epsilon\}
\end{equation*}
is relatively dense in $\R$, i.e. there exists a number $l>0$ such
that
\begin{equation*}
T(g, \epsilon)\bigcap [t , t+l]\neq\emptyset, \quad \forall t\in \R.
\end{equation*}
As is well known, for $g$ to be almost periodic it is necessary and
sufficient that the family of functions $\{S(t)g\}_{t\in \R}$ is
pre-compact in $BUC(\R, \X)$.

For $g\in AP(X)$, let
\begin{equation*}
a(\lambda, g) := \lim_{T\rightarrow
\infty}\frac{1}{2T}\int_{-T}^{T}e^{-i\lambda t}g(t)dt, \quad
\forall \lambda \in \R.
\end{equation*}
It is known that $ a(\lambda, g) $ is well defined and there are at
most countably many points $\lambda$ such that $a(\lambda, g)\neq
0$. We call the set $\sigma_b(g) := \{\lambda :a(\lambda, g) \neq
0\}$ and $\mathfrak{M} := \{\sum_{k=1}^{N}n_k\lambda_k : \lambda_k
\in \sigma_b(g), n_k\in \Z, \forall N\in \N\}$ the Bohr spectrum and
module of $g$, respectively.

Under the above notations, $g$ can be approximated uniformly on $\R$
by a sequence of trigonometric polynomial(see\cite[ Chap.
2]{levzhi})
\begin{equation*}
P_n(t) := \sum_{k=1}^{N(n)}a_{n,k}e^{i\lambda_{n,k}t},\quad n=1, 2,
\cdots; \quad \lambda_{n,k} \in \sigma_b(g), \quad t\in\R.
\end{equation*}
And of course every function which can be approximated uniformly by
a sequence of trigonometric polynomial is almost periodic.

A finite or countable set of real numbers $\Omega := \{\omega_1,
\omega_2, \cdots, \omega_k, \cdots\}$ is said to be rationally
independent if for any $n\in \N$, the following relation holds
\begin{equation*}
\left.\begin{matrix}
 r_1\omega_1 + r_2\omega_2 + \cdots + r_n\omega_n = 0\\
 r_1, r_2, \cdots, r_n\in \mathbb{Q}
\end{matrix}
\right\}
 \Rightarrow r_1 = r_2 = \cdots
= r_n =0.
\end{equation*}
We call $\Omega$ a rational basis of $\sigma_b(g)$ if it is
rationally independent and every $\lambda_j\in \sigma_b(g)$ is
representable as a finite linear combination of the $\omega_k$ with
rational coefficients, that is,
\begin{equation*}
\lambda_j = r^{(j)}_1\omega_1 + r^{(j)}_2\omega_2 + \cdots +
r^{(j)}_n\omega_n, \quad r^{(j)}_1, r^{(j)}_2, \cdots,
r^{(j)}_{n_j}\in\mathbb{Q}, \quad (j = 1, 2, \cdots).
\end{equation*}
Moreover, if all the $r^{(j)}_k\in \Z$ , then the basis is called
integer basis.

\begin{definition}
If $\mathfrak{M}(g)$ has a finite integer basis $\{\omega_1,
\omega_2, \cdots, \omega_k\}$, then $g$ is called a quasi-periodic
function with frequencies $(\omega_1/2\pi, \omega_2/2\pi, \cdots,
\omega_k/2\pi)$. In the sequel, we also call $g$ {\it
k-basic-frequency-quasi-periodic function}, or {\it
$k$-quasi-periodic function} .
\end{definition}
From the definition it follows that if $g$ is $k$-quasi-periodic, it
cannot be $l$-quasi-periodic for any $l<k$.

\subsection{Spectral Theory of Functions}
 Let $g\in BC(\R,\X)$, and let $\lambda \in
\C$ such that $Re\lambda \not= 0$. Then, obviously, the equation $
\dot x = \lambda x , \quad x\in \R $ has an exponential dichotomy.
By the Perron Theorem in ODE, its non-homogeneous equation $\dot
x=\lambda x+g(t)$ has a unique solution $x_{g,\lambda}\in
BUC(\R,\X)$. Moreover, $x_{g,\lambda}= ({\cal D}-\lambda )^{-1}g$
for every $g\in BC(\R,\X)$. Therefore, ${\cal D}$ is a closed
operator on $BC(\R,\X)$ and $\rho ({\cal D}) \supset
\mathbb{C}\backslash i\R$, and $(\lambda -{\cal D})^{-1}g$, as a
function of $\lambda$, is analytic everywhere in $\C\backslash i\R$.
\begin{definition}
The set of all reals $\xi \in \R$ such that the complex function
$(\lambda -{\cal D})^{-1}g$ has no analytic extension to any
neighborhood of $i\xi$ is said to be the uniform spectrum of $g$,
and is denoted by $sp_u(g)$.
\end{definition}
Using the Green function to determine bounded solutions in the
theory of ODE we have
\begin{eqnarray}\label{re-2}
(\lambda -{\cal D})^{-1}g (\xi )&=& -x_{g ,\lambda }(\xi ) \nonumber\\
&=&
\begin{cases}
\begin{array}{ll}
\int^  {\infty}_\xi  e^{\lambda (\xi -t)}g(t)dt,     &(\mbox{if}\
Re\lambda > 0),\\ \\
-\int _ {-\infty}^\xi  e^{\lambda (\xi -t)}g(t)dt,    & (\mbox{if }
\ Re\lambda < 0)
\end{array} \end{cases} \nonumber\\
&=&  \begin{cases}
\begin{array}{ll}
\int^  {\infty}_0  e^{-\lambda \eta }g(\xi +\eta )d\eta,
&(\mbox{if}\ Re\lambda > 0),\\ \\
-\int _ {-\infty}^0  e^{-\lambda \eta}g(\xi +\eta )d\eta ,    &
(\mbox{if } \ Re\lambda < 0).
\end{array} \end{cases}
\end{eqnarray}
By definition,
\begin{eqnarray}
\hat{g} (\lambda ) := (\lambda -{\cal D})^{-1}g (0 )&=&
\begin{cases}
\begin{array}{ll} \int^  {\infty}_0  e^{-\lambda \eta }g(\eta
)d\eta,
&(\mbox{if}\ Re\lambda > 0),\\ \\
-\int _ {-\infty}^0  e^{-\lambda \eta}g(\eta )d\eta,     & (\mbox{if
} \ Re\lambda < 0)
\end{array} \end{cases}
\end{eqnarray}
is called the Carleman-Laplace transform of $g$. Obviously, the
Carleman-Laplace transform of $g$ is analytic in $\lambda \in
\C\backslash i\R$.
\begin{definition}
Let $g\in BC(\R,\X)$. The set of all reals $\xi$ such that
$\hat{g}(\lambda )$ has no analytic extension to any neighborhood of
$i\xi$ is called the Carleman spectrum of $g$ and is denoted by
$sp_c(g)$.
\end{definition}
It is easy to see from the definitions that $sp_c(g)\subset
sp_u(g)$. In \cite{liunguminvu}, it is proved that they are actually
coincide. As is well known (see e.g. \cite{pru}), the Carleman
spectrum $sp_c(g)$ coincides with the Beurling spectrum $sp_b(g)$
defined as
\begin{equation}
sp_b(g):= \{ \xi\in\R :\ \forall \epsilon >0, \exists \varphi\in
L^1(\R), \ supp \tilde{\varphi} \subset (\xi -\epsilon ,\xi+\epsilon
), \ \varphi * g \not= 0 \} ,
\end{equation}
where $\tilde{\varphi}$ is the Fourier transform of $\varphi$,
that is,
$$
\tilde{\varphi}(\eta )= \int^\infty_{-\infty} e^{-i\eta t}\varphi
(t)dt ,\quad \eta\in\R .
$$
For this reason, for each $g\in BC(\R,\X)$ we will denote simply by
$sp(g)$ the spectrum of $g$ in any sense mentioned above.

\bigskip
Below we give some properties of the spectrum of functions which we
will need in the sequel.

\begin{proposition}
Let $g, g_n\in BUC(\R, \X), n\in \mathbb{N}$, such that
$g_n\rightarrow g$ as $n\rightarrow\infty$. Then
\begin{enumerate}
\item $sp(g)$ is closed,

\item $sp(g(\cdot+h)) = sp(g)$,

\item If $\alpha \in \mathbb{C}\backslash\{0\}$, then $sp(\alpha
g)=sp(g)$,

\item If $sp(g_n)\subset\Lambda$ for all $n\in \N$, then
 $sp(g)\subset\bar{\Lambda}$,

\item If $A$ is a closed operator , $g(t)\in D(A), \forall t\in
\R$ and $\mathcal{A}g \in BUC(\R, \X)$, then
$sp(\mathcal{A}g)\subset
  sp(g)$.
  \end{enumerate}
\end{proposition}
As a consequence of this proposition, for a closed subset of $\R$,
say, $\Lambda$, then the following
\begin{eqnarray*}
\Lambda _{BUC}(\X) &:=& \{ g\in BUC(\R,\X):\ sp(g)\subset \Lambda\},\\
\Lambda _{AP}(\X) &:=& \{ g\in AP(\X):\ sp(g)\subset \Lambda \} ,
\end{eqnarray*}
are closed translation-invariant subspaces of $BC(\R,\X)$.
Considering the translation group $(S(t))_{t\in \R}$ on $BUC(\R,
\X)$, we give a frequently used properties of spectrum.
\begin{lemma}\label{lem 11}
Denote ${\cal D}_\Lambda$ the restriction of $\mathcal{D}$ on
$\Lambda_{BUC}(\X)$. Then
\begin{equation*}
\sigma({\cal D}_\Lambda)=i\Lambda.
\end{equation*}
\end{lemma}
\begin{proof}
See  \cite[Lemma 3.3]{murnaimin}.
\end{proof}

\bigskip
We also recall that the reduced spectrum $sp_{AP}(g)$ is defined as
follows:
\begin{equation*}
sp_{AP}(g):= \{ \xi\in\R :\ \forall \epsilon >0, \exists \varphi\in
L^1(\R), \ supp \tilde{\varphi} \subset (\xi -\epsilon ,\xi+\epsilon
), \ \varphi * f \not\in AP(\X) \} .
\end{equation*}
Let us consider the quotient space $BUC(\R, \X)/AP(\X) $. Since
$(S(t))_{t\in\R}$ leaves $AP(\X)$ invariant there is an induced
translation group $(\tilde{S}(t))_{t\in\R}$ on $BUC(\R, \X)/AP(\X)
$.
 $(\tilde{S}(t))_{t\in\R}$ and its generator $\tilde{\mathcal{D}}$ are given by
\begin{eqnarray*}
\tilde{S}(t)\pi(g) &=& \pi(S(t)g),\quad t\in\R,g\in BUC(\R, \X),\\
\tilde{\mathcal{D}}\pi(g) &=& \pi(\mathcal{D}g),\qquad g\in
D(\mathcal{D}),
\end{eqnarray*}
where $\pi : BUC(\R, \X)\rightarrow BUC(\R, \X)/AP(\X)$ denotes the
quotient mapping. For $g\in BUC(\R,\X)$,  as in \cite[the proof of
Theorem 4.3, p. 374]{arebat}, $ sp_{AP}(g)$ coincides with the set
of all reals $\xi$ such that $(\lambda -
\tilde{\mathcal{D}})^{-1}\pi(g)$ has no analytic extension to any
neighborhood of $i\xi$. The following is the Loomis Theorem in the
infinite dimensional case.

\begin{theorem}\label{th 11}
Let $\X$ be a Banach space that does not contain any subspaces
isomorphic to $c_0$, and let $g$ be in $ BUC(\R,\X)$ with
countable $sp_{AP}(g)$. Then, $g\in AP(\X)$.
\end{theorem}
\begin{proof}
For the proof we refer the reader to \cite{arebat,arebathieneu}.
\end{proof}


\subsection{Decomposition of Function Spaces}
In this subsection we will decompose a function space into direct
sum of simpler ones in terms of spectrum.

\begin{lemma}\label{lem 1}
Suppose the closed subset $\Lambda\subset\R$ can be represented as
the disjoint union of a compact subsets $\Lambda^1$ and a closed
set $\Lambda^2$, i.e.,
\begin{equation*}
\Lambda= \Lambda^1 \cup\Lambda^2 .
\end{equation*}
Then,
\begin{equation*}
\Lambda_{BUC}(\X) =\Lambda^1_{BUC}(\X)\oplus \Lambda^2_{BUC}(\X) .
\end{equation*}
Moreover, the projections $P,Q$ corresponding to this splitting
are determined from $R(\lambda, {\cal D})$ that commute with any
operators that commute with ${\cal D}$ in the sense of resolvents
commuting.
\end{lemma}
\begin{proof}
By Lemma \ref{lem 11}, we have
\begin{equation*}
\sigma({\cal D}_\Lambda)=i\Lambda=i\Lambda^1 \cup i\Lambda^2.
\end{equation*}
In case ${\cal D}_\Lambda$ is bounded we can prove the lemma by
using the Riesz projection. For example, the subset $\Lambda$ is
compact. Now, suppose ${\cal D}_\Lambda$ is unbounded. Then, for a
fixed $\mu\in\rho({\cal D}_\Lambda)$, we have $0\in \sigma(R(\mu,
{\cal D}_\Lambda))$ and $\sigma(R(\mu, {\cal D}_\Lambda))\setminus
\{0\}=(\mu-\sigma({\cal
D}_\Lambda))^{-1}:=\{\frac{1}{\mu-\lambda}:\lambda\in \sigma({\cal
D}_\Lambda)\}$ (see \cite[Chapter iv, sect. 1, p. 227]{engnag}).
Therefore, we obtain
\begin{align*}
\sigma(R(\mu, {\cal D}_\Lambda))&=(\mu-\sigma({\cal
D}_{\Lambda^1}))^{-1}\cup(\mu-\sigma({\cal
D}_{\Lambda^2}))^{-1}\cup \{0\}\\
 &=:\tau_1\cup\tau_2,
\end{align*}
where $\tau_1:=(\mu-\sigma({\cal D}_{\Lambda^1}))^{-1},
\tau_2:=(\mu-\sigma({\cal D}_{\Lambda^2}))^{-1}\cup \{0\}$ are
compact and disjoint subsets of $\mathbb{C}$.

Now let
\begin{equation*}
P:=\frac{1}{2\pi i}\int_{\gamma}R(\lambda,R(\mu, {\cal
D}_\Lambda))d\lambda,
\end{equation*}
where $\gamma$ is a Jordan path in the complement of $ \tau_2$ and
enclosing $\tau_1$. This projection commutes with any operators
that commute with ${\cal D_\Lambda}$ in the sense of resolvents
commuting and yields the spectral decomposition
\begin{equation*}
\Lambda_{BUC}(\X)=\X_1\oplus \X_2
\end{equation*}
with ${\cal D_\Lambda}$-invariant spaces $\X_1=P\Lambda_{BUC}(\X),
X_2=Q\Lambda_{BUC}(\X),Q:=I-P.$ The restrictions ${\cal D}_{\X_1}$
and ${\cal D}_{\X_2}$ of ${\cal D}_{\Lambda}$ satisfy
\begin{equation*}
\sigma({\cal D}_{\X_1})=i\Lambda^1, \quad \sigma({\cal
D}_{\X_2})=i\Lambda^2.
\end{equation*}
Let $\xi\in \R\backslash \Lambda^1.$ Then, for every $  g\in \X_1,
(\lambda-{\cal D})^{-1}g=(\lambda-{\cal D}_{\X_1})^{-1}g$ has an
analytic extension to some neighborhood of i$\xi$. So
$sp(g)\subset\Lambda^1$, and $\X_1\subset\Lambda^1_{BUC}(\X)$.
Similarly, $\X_2\subset\Lambda^2_{BUC}(\X)$. On the other hand,
since $\Lambda^1_{BUC}(\X)\bigcap\Lambda^2_{BUC}(\X)=\{0\}$, we
obtain
\begin{equation*}
\Lambda_{BUC}(\X) =\Lambda^1_{BUC}(\X)\oplus \Lambda^2_{BUC}(\X) .
\end{equation*}
\end{proof}
When both sets $\Lambda^1$ and $\Lambda^2$ are not compact we have
the following decomposition. Let $\Sigma,\Sigma_1,\Sigma_2$ be
closed subsets of the unit circle $\Gamma$ such that $\Sigma$ is
the disjoint union of $ \Sigma_1$ and $\Sigma_2$ . Denote
\begin{eqnarray}
\Lambda  &=& \{\xi\in\R:\ e^{i\xi}\in \Sigma \},\\
\Lambda^1 &=& \{\xi\in\R:\ e^{i\xi}\in \Sigma_1 \},\\
\Lambda^2 &=& \{\xi\in\R:\ e^{i\xi}\in \Sigma _2\} .
\end{eqnarray}
\begin{lemma}\label{lem 7}
Under the above notations, one has
\begin{equation}
\Lambda_{BUC}(\X) = \Lambda^1 _{BUC}(\X) \oplus
\Lambda^2_{BUC}(\X) .
\end{equation}
Moreover, the projections $P,Q$ corresponding to this splitting are
determined as the Riesz spectral projections of the translation
$S(1)$ on $\Lambda_{BUC}(\X)$ with respect to $\Sigma_1$ and
$\Sigma_2$, so each operator commuting with translation commutes
with $P,Q$ as well.
\end{lemma}
\begin{proof}
Consider the operators $\mathcal{D}_\Lambda$ on
$\Lambda_{BUC}(\X)$, $(S_\Lambda(t))_{t \in \R}$ on
$\Lambda_{BUC}(\X)$. By the Weak Spectral Mapping Theorem, we
obtain
\begin{equation*}
\sigma(S_\Lambda(1))) =
\overline{e^{\sigma(\mathcal{D}_{\Lambda_{BUC}(\X)})}} =
\overline{e^{i\Lambda}} =
 \Sigma = \Sigma_1\cup\Sigma_2.
\end{equation*}
Using the Riesz projection $P$ on $\Lambda_{BUC}(\X)$ for $S(1)$, we
have
\begin{equation*}
\Lambda_{BUC}(\X) = P\Lambda_{BUC}(\X) \oplus Q\Lambda_{BUC}(\X),
\end{equation*}
where $P := \frac{1}{2\pi i}\int_\gamma R(\lambda,
S_\Lambda(1))d\lambda$, $Q := I - P$, $\gamma$ is a simple contour
in the complement of $ \Sigma_2$ and enclosing $\Sigma_1$. We will
prove that $P\Lambda_{BUC}(\X) = \Lambda^1_{BUC}(\X)$. In fact, by
the Weak Spectral Mapping Theorem
\begin{equation*}
\overline{e^{\sigma(\mathcal{D}|_{P\Lambda_{BUC}(\X)})}} =
\sigma(S_\Lambda(1))|_{P\Lambda_{BUC}(\X)})  = \Sigma_1 =
\overline{e^{i\Lambda^1}}.
\end{equation*}
So
\begin{equation*}
\sigma(\mathcal{D}|_{P\Lambda_{BUC}(\X)}) \subset i\Lambda^1.
\end{equation*}
Therefore, for each $u\in P\Lambda_{BUC}(\X)$ one has $sp(u)
\subset \Lambda^1.$ Similarly, if $w\in Q\Lambda_{BUC}(\X)$, then
$sp(w) \subset \Lambda^2.$ Obviously,
$\Lambda^1_{BUC}(\X)\bigcap\Lambda^2_{BUC}(\X) = \{0\}$ and they
are closed subspaces of $\Lambda_{BUC}(\X)$. So, we have
\begin{equation*}
\Lambda_{BUC}(\X) = \Lambda^1_{BUC}(\X) \oplus \Lambda^2_{BUC}(\X)
\end{equation*}
and the corresponding projection $P = \frac{1}{2\pi i}\int_\gamma
R(\lambda, S_\Lambda(1))d\lambda$. Note that $P$ is determined from
$S(t)$, so it commutes with any operators that commute with $S(t)$,
or $\mathcal{D}$.
\end{proof}

\section{Main Results}

We begin this section by recalling some concept of solutions.
\begin{definition}
A function $u\in BUC(\R,\X)$ is said to be a classical solution of
Eq.(\ref{1}) if $\dot u(t)$ exists as an element of $BUC(\R,\X)$
such that for all $t$, $u(t) \in D(A)$, and Eq.(\ref{1}) holds.
\end{definition}

\begin{definition}\label{def 2}
A function $u\in BUC(\R,\X)$ is said to be a mild solution of
Eq.(\ref{1}) if for all $t\in \R$, $\int^t_0 u(\xi )d\xi \in D(A)$,
and the following equation is satisfied:
\begin{equation}\label{def mild solution}
u(t)-u(0) = A\int^t_0 u(\xi )d\xi + \int^t_0 [{\cal B}u(\xi ) +f(\xi
)]d\xi .
\end{equation}
\end{definition}
Obviously, every classical solution is a mild solution, but not vice
versa in general. We can check that a mild solution $u$ that is
continuously differentiable, and $u(t)\in D(A)$ for all $t$, is a
classical one. We note that (see e.g.
\cite{arebathieneu,hinnaiminshi,murnaimin}) if $A$ generates a
semigroup $(T(t))_{t\ge 0}$, then it is a mild solution in the sense
that it satisfies the following integral equation
$$
u(t) = T(t-s) u(s) + \int^t_s T(t-\xi ) [{\cal B}u(\xi ) + f(\xi
)] d\xi , \quad \forall t\ge s .
$$

In the sequel, if $A$ is an operator in a Banach space $\Y$, we
use the notation
$$
\sigma_i (A):= \{ \xi \in \R :\ i\xi \in \sigma (A)\} .
$$
\begin{theorem}\label{th 3}
Let Eq.(\ref{1}) have a mild solution $u$ in $BUC(\R,\X)$,
 $\sigma _i ({\cal A}+{\cal B})\backslash sp(f)$ be compact, $sp(f)$ be countable, and  $\X$ do not contain any subspace isomorphic
 to $c_0$. Then Eq.(\ref{1}) has an almost periodic mild solution
 $u_f$ with $sp(u_f) \subset sp(f)$.
\end{theorem}
We would like to present some useful results before giving the
proof of Theorem \ref{th 3}.

\medskip
\begin{definition}
Let $\Lambda$ be any closed subset of the real line. We define an
operator $L_\Lambda$ as follows: $u \in D(L_\Lambda )$ if and only
if $u\in \Lambda_{BUC}(\X), \int_{0}^{t}u(\xi)d\xi \in D(A)$ for all
t and there is $g\in \Lambda_{BUC}(\X)$ such that
$$
u(t)-u(0) = A\int^t_0 u(\xi )d\xi + \int^t_0 g(\xi )d\xi ,\quad
\forall  t\in \R ,
$$
and in this case $L_\Lambda u := g$.
\end{definition}

At this point, we can check easily that $L$ is a linear operator,
but we do not know if this operator is a single-valued operator. The
next lemma will clarify this. Note that since ${\cal B}$ is bounded
and commutes with translations, one can easily show that it leaves
$\Lambda_{BUC}(\X)$ invariant. Therefore, by abuse of notation we
will still use ${\cal B}$ to denote its restriction to
$\Lambda_{BUC}(\X)$ if this does not cause any danger of confusion.

\begin{lemma}\label{lem 3}
The operator $L_\Lambda $ is a linear single-valued and closed
operator on $\Lambda_{BUC}(\X)$, so $L_\Lambda -{\cal B}$ is a
closed single-valued operator in $\Lambda_{BUC}(\X)$.
\end{lemma}
\begin{proof}
First, we show that $L_\Lambda$ is a single-valued closed operator.
Suppose $u\in D(L_\Lambda), g_1, g_2\in \Lambda_{BUC}(\X)$ such that
$$
L_\Lambda u = g_1, \quad L_\Lambda u = g_2,
$$
then
$$
\int^t_0 (g_1(\xi )-g_2(\xi))d\xi = 0 ,\quad \forall t ,
$$
thus
$$
\frac{1}{t-s}\int^t_0 (g_1(\xi )-g_2(\xi))d\xi = 0 ,\quad \forall
t>s .
$$
Fixing $t$ and letting $s\rightarrow t$ we obtain $g_1(t
)-g_2(t)=0$, i.e. $L_\Lambda$ is single-valued.

Suppose $u, u_n \in D(L_\Lambda ), \forall n \in \N, g, g_n \in
\Lambda_{BUC}(\X), \forall n \in \N, $ such that $L_\Lambda u_n=g_n$
and
\begin{equation*}
\lim_{n\rightarrow\infty} u_n = u,\quad \lim_{n\rightarrow\infty}
g_n =g
\end{equation*}
 in the sup-norm topology of $BUC(\R, \X) $, then
$$
\int^{t}_{0}u_n(\xi)d\xi \rightarrow \int^{t}_{0}u(\xi)d\xi,\qquad
(n\rightarrow\infty) \quad \forall t,
$$
$$
\int^{t}_{0}g_n(\xi)d\xi \rightarrow \int^{t}_{0}g(\xi)d\xi,\qquad
(n\rightarrow\infty) \quad \forall t.
$$
By the closedness of $A$, for each fixed $t\in\R$ we have
$$
u(t)-u(0) = A\int^t_0 u(\xi )d\xi + \int^t_0 g(\xi )d\xi ,\quad
\forall  t .
$$
This shows that $L_\Lambda$ is closed.

Finally, since ${\cal B}$ is bounded is bounded $L_\Lambda -{\cal
B}$ is a closed operator on $\Lambda_{BUC}(\X )$. \end{proof}

For the relation among the spectrum of mild solutions,
$\mathcal{A+B}$ and $f$ we have the following theorem that has been
known for abstract ODE (see e.g. \cite{arebathieneu,levzhi}).
\begin{theorem}\label{th 1}
Let $u$ be a mild solution of Eq.(\ref{1}) on $\R$. Then
\begin{equation*}
sp(u) \subset \sigma _i ({\cal A}+{\cal B}) \cup sp(f),
\end{equation*}
where $\sigma _i ({\cal A}+{\cal B}) := \{\xi\in \R: i\xi \in \sigma
({\cal A}+{\cal B})\}$.
\end{theorem}
\begin{proof}
For $Re \lambda >0$ taking the Carleman-Laplace transform of both
sides of Eq. (\ref{def mild solution}) we have
\begin{equation}
\hat u(\lambda ) -\frac{1}{\lambda}u(0) = \frac{1}{\lambda} A \hat
u(\lambda ) +\frac{1}{\lambda}  \hat{[{\cal B} u]}(\lambda )
+\frac{1}{\lambda}  \hat f(\lambda ).
\end{equation}
Therefore,
\begin{equation}
\lambda \hat u(\lambda ) -u(0) =  A \hat u(\lambda ) + \hat{[{\cal
B} u]}(\lambda ) + \hat f(\lambda ).
\end{equation}
Similarly, since $S(s)u$ is again a mild solution of Eq.(\ref{1})
with $f$ replaced by $S(s)f$ for each fixed $s\in\R$, we have
\begin{equation}
\lambda \widehat{S(s) u}(\lambda ) -[S(s)u](0) =  A \widehat{S(s)
u}(\lambda ) + \widehat{[{\cal B} S(s)u]}(\lambda ) + \widehat{S(s)
f}(\lambda ).
\end{equation}
By (\ref{re-2}), this means
\begin{equation}
\lambda R(\lambda ,{\cal D})u - u =  {\cal A}R(\lambda ,{\cal D})u +
R(\lambda ,{\cal D}){\cal B}u + R(\lambda ,{\cal D})f.
\end{equation}
Since ${\cal B}$ commutes with ${\cal D}$, we have
\begin{equation}\label{3.6}
\lambda R(\lambda ,{\cal D})u - u =  {\cal A}R(\lambda ,{\cal D})u +
{\cal B}R(\lambda ,{\cal D}) u + R(\lambda ,{\cal D})f.
\end{equation}

Suppose $\xi$ is a real such that $\xi \not\in sp(f)$, and $i\xi
\not\in \sigma ({\cal A}+{\cal B})$. Then $R(\lambda ,{\cal D})f$,
and $\left[ \lambda -{\cal A} - {\cal B}\right]^{-1}$ are extendable
to analytic functions around $i\xi$. And for such a $\xi$ we have
\begin{eqnarray*}
(\lambda  -{\cal A} - {\cal B})^{-1}\left[ R(\lambda ,{\cal D})f
+u\right]  &=& R(\lambda ,{\cal D}) u .
\end{eqnarray*}
Therefore, for such a $\xi$, $R(\lambda ,{\cal D})u$ has an analytic
extension to a neighborhood of $i\xi$. This show that $ \xi \not\in
sp(u)$. The theorem is proved.
\end{proof}
Lemma \ref{lem 3} and Theorem \ref{th 1} yield the following

\begin{lemma}\label{lem 2}
Let  $\Lambda =\sigma _i ({\cal A}+{\cal B})\cup sp(f)$. Then a
function $u\in BUC(\R ,\X)$ is a mild solution of Eq.(\ref{1}) if
and only if
\begin{equation}\label{cl}
(L_\Lambda -{\cal B}) u=f.
\end{equation}
\end{lemma}
\begin{proof}
By Lemma \ref{lem 3} the operator $L_\Lambda$ is well-defined as a
single-valued operator. Obviously, if $u$ is a solution of Eq.
(\ref{cl}), then by the definition of the operator $L_\Lambda$ it
satisfies Eq. (\ref{def mild solution}), so it is a mild solution of
Eq. (\ref{1}).

\medskip
Conversely, let $u\in BUC(\R ,\X)$ be a mild solution of Eq.
(\ref{1}). Then, by Theorem \ref{th 1} it is in $\Lambda_{BUC}(\X)$,
so, by the definition of the operator $L_\Lambda$, it satisfies Eq.
(\ref{cl}).
\end{proof}

Noticing that $f$ is almost periodic, similarly as the proof of
Theorem \ref{th 1}, we can obtain the following theorem about the
reduced spectrum of mild solutions.

\begin{theorem}\label{th 2}
Let $u$ be a mild solution of Eq.(\ref{1}) on $\R$. Then,
\begin{equation}
sp_{AP}(u) \subset \sigma _i ({\cal A}+{\cal B}).
\end{equation}
\end{theorem}
\begin{proof}
 Recall that ${\cal A}+{\cal B}$ is an operator defined in
$BUC(\R,\X)$. By (\ref{3.6}),
\begin{equation}\label{3.9}
[\lambda  - {\cal A} - {\cal B}]R(\lambda ,{\cal D}) u = R(\lambda
,{\cal D})f +u.
\end{equation}
Let $\xi\in \R$ such that $i\xi \in \rho ({\cal A}+{\cal B})$. Then,
the resolvent $R(\lambda , {\cal A}+{\cal B})$ exists, and is
analytic in a neighborhood of $i\xi$ in $\C$. Therefore, in a
neighborhood of $i\xi$ if $Re\lambda \not=0$, by (\ref{3.9}) we have
\begin{equation}\label{3.10}
 R(\lambda ,{\cal D}) u = R(\lambda , {\cal A}+{\cal B}) (R(\lambda
,{\cal D})f +u).
\end{equation}
Note that each bounded linear operator in $BUC(\R,\X)$ that commutes
with the translations must map $AP(\X)$ into itself. Therefore, for
$Re\lambda\not= 0$, $ R(\lambda ,{\cal D})$ maps $AP(\X)$ into
itself. From the autonomousness of ${\cal B}$, it follows also that
$R(\lambda , {\cal A}+{\cal B})$ commutes with the translations, so
it leaves $AP(\X)$ invariant. This means that the operators $
R(\lambda ,{\cal D})$ and $R(\lambda , {\cal A}+{\cal B})$
 induce operators $\widetilde{R(\lambda ,{\cal D})}= R(\lambda ,\tilde{{\cal D}})$ and
$\widetilde{R( \lambda ,{\cal A} + {\cal B})}$ in
$BUC(\R,\X)/AP(\X)$, respectively. Since $f\in AP(\R,\X)$, we have $
R(\lambda ,\tilde{{\cal D}}) \pi(f)=0.$ And hence, by (\ref{3.10})
we have
\begin{eqnarray*}
R(\lambda ,\tilde{{\cal D}})\pi(u) = \widetilde{R( \lambda ,{\cal A}
+ {\cal B})}  \pi(u) .
\end{eqnarray*}
Therefore, $\widetilde{R( \lambda ,{\cal A} + {\cal B})}  \pi(u)$ is
a natural analytic extension to a neighborhood of $i\xi$ for $
R(\lambda ,\tilde{{\cal D}})\pi(u)$. This shows that $\xi \not\in
sp_{AP}(u)$ and the proof is completed.
\end{proof}

By Theorem \ref{th 11} and Theorem \ref{th 2}  we obtain the
following corollary immediately.

\begin{corollary}
Let
 $\sigma _i ({\cal A}+{\cal B})$ be countable. Then every bounded
 and uniformly continuous mild solution of Eq.(\ref{1}) is almost
 periodic provided $\X$ does not contain any subspace isomorphic
 to $c_0$.
\end{corollary}

The following lemma is standard (see \cite[Proposition
III.6.5]{kat}).
\begin{lemma}\label{lem 4}
Let $A$ be a closed operator on a Banach space $\X$ such that
$\rho (A)\not= \emptyset$, and let $B$ be a bounded operator on
$\X$. Then, $A$ commutes with $B$ in the sense of resolvents
commuting if and only if they commute in the sense that $BD(A)
\subset D(A)$ and $ABx=BAx$ for all $x\in D(A)$.
\end{lemma}

Now we are ready to complete the proof of Theorem \ref{th 3}.

\noindent{\it Proof of Theorem \ref{th 3}.} Let $\Lambda^1:=\sigma
_i ({\cal A}+{\cal B})\backslash sp(f)$ and $\Lambda^2:= sp(f)$. By
Theorem \ref{th 1}, $sp(u)\subset \Lambda := \Lambda^1\cup
\Lambda^2$. By the compactness of $\Lambda^1$ and Lemma \ref{lem 1}
we have
$$
\Lambda_{BUC}(\X) =\Lambda^1_{BUC}(\X)\oplus \Lambda^2_{BUC}(\X),
$$
and the projections $P,Q$ corresponding to this splitting commute
with each operator that commutes with ${\cal D}_\Lambda$ in the
sense of resolvents commuting. Note that $L_\Lambda$ commutes with
all translations, so it commutes with ${\cal D}_\Lambda$. By Lemma
\ref{lem 2} we have
$$
Q[L_\Lambda -{\cal B}] u=Qf=f.
$$
By Lemma \ref{lem 4}, $Q[L_\Lambda -{\cal B}]u=[L_\Lambda -{\cal B}]
Qu$, and letting $Qu=w$, we have
$$
[L_\Lambda -{\cal B}] w= f.
$$
That is, $w\in \Lambda^2_{BUC}(\X)$ is a mild solution of
Eq.(\ref{1}). Since $\Lambda^2$ is countable and $\X$ does not
contain any subspace isomorphic to $c_0$, by Theorem \ref{th 11},
$w$ is almost periodic. The theorem is proved. \qed

\bigskip
Now we consider the case when $\sigma _i ({\cal A}+{\cal
B})\backslash sp(f)$ may not be compact.
\begin{theorem}\label{th 7}
Suppose that Eq.(\ref{1}) has a mild solution $u\in BUC(\R, \X)$,
and $\overline{e^{i\sigma_i(\mathcal{A} + B)}}\backslash
\overline{e^{isp(f)}}$ is closed. Then, there exists a mild solution
$w$ of Eq.(\ref{1}) such that
\begin{equation*}
e^{isp(w)}\subset \overline{e^{isp(f)}} .
\end{equation*}
\end{theorem}
\begin{proof}
Let $\Sigma := \overline{e^{i\sigma_i(\mathcal{A} + B)}}\cup
\overline{e^{isp(f)}}, \Sigma_1 :=
\overline{e^{i\sigma_i(\mathcal{A} + B)}}\backslash
\overline{e^{isp(f)}}, \Sigma_2 := \overline{e^{isp(f)}}$.
Obviously, $\Sigma,\Sigma_1,\Sigma_2$ are closed. Moreover, by the
assumption, $\Sigma_1\cap\Sigma_2=\emptyset $. Therefore, by Lemma
\ref{lem 7}, if we define
\begin{eqnarray*}
\Lambda  &=& \{\xi\in\R:\ e^{i\xi}\in \Sigma \},\\
\Lambda^1 &=& \{\xi\in\R:\ e^{i\xi}\in \Sigma_1 \},\\
\Lambda^2 &=& \{\xi\in\R:\ e^{i\xi}\in \Sigma _2\} ,
\end{eqnarray*}
then
\begin{equation*} \Lambda_{BUC}(\X) = \Lambda^1_{BUC}(\X)
\oplus \Lambda^2_{BUC}(\X) .
\end{equation*}
Moreover, the projections $P$ and $Q$ corresponding to this
splitting commute with any operators that commute with translations.
In the same way as in the proof of Theorem \ref{th 3} we can show
that $w:=Qu$ is the sought mild solution.
\end{proof}

Now we prove a Massera type theorem for almost periodic solutions
of Eq.(\ref{1}) when $\sigma _i ({\cal A}+{\cal B})\backslash
sp(f)$ may not be compact.
\begin{theorem}\label{th 8}
In addition to the conditions of Theorem \ref{th 7} assume that
$\X$ does not contain any subspace isomorphic
 to $c_0$ and $\overline{e^{isp(f)}}$ is countable. Then Eq.(\ref{1}) has an almost periodic mild solution $w$ such that
 $$
e^{isp(w)} \subset \overline{e^{isp(f)}}.
 $$
\end{theorem}
\begin{proof}
By Theorem \ref{th 7}, we obtain that $e^{isp(w)}$ is countable.
This implies that $sp(w)$ is countable. By Loomis'  Theorem we have
$w\in AP(\X)$.
\end{proof}

\bigskip
Let Eq. (\ref{1}) have an almost periodic mild solution. We will
refine the technique of decomposition to further the results
obtained above. We first give some lemmas. The following is an easy
lemma whose proof is given for the reader's convenience.
\begin{lemma}
Let $g\in AP(\X)$ and let $\cal B$ be an autonomous operator in
$BUC(\R,\X)$. Then $\mathcal{B}g\in AP(\X)$.
\end{lemma}
\begin{proof}
We will use the Bochner criterion for the almost periodicity of a
function. Since $g\in AP(\X)$, we have that the set $\{S(t)g, t\in
R\}$ is pre-compact in $BUC(\R,\X)$. On the other hand, as ${\cal
B}$ is a bounded linear {\it autonomous} operator in $BUC(\R,\X)$,
we obtain that the set $ \{S(t)\mathcal{B}g, t\in
R\}=\mathcal{B}\{S(t)g, t\in R\} $ is pre-compact in $BUC(\R,\X)$,
that is $\mathcal{B}g\in AP(\X)$.
\end{proof}

\begin{lemma}\label{lem 5}
Denote $\Y := BUC(\R,\X)$. For any $g\in \Y$, we define the function
$G$ as follows:
\begin{eqnarray*}
G &:& \R \rightarrow \Y\\
G &:& t\mapsto S(t)g.
\end{eqnarray*}
If $g \in AP(\X)$, then $G \in AP(\Y)$  and
$\sigma_{b}(G)=\sigma_{b}(g)$.
\end{lemma}
\begin{proof}
Since $g$ is almost periodic, $T(g, \epsilon)$ is relatively dense
in $\R$ for all $\epsilon > 0$. On the other hand,
\begin{eqnarray*}
\sup_{t \in \R}|G(t+\tau)-G(t)| = \sup_{t \in \R}|S(t+\tau)g-S(t)g|=
\sup_{s\in \R}\|g(\tau + s)-g(s)\|,
\end{eqnarray*}
 where  $\|\cdot \|$ and $|\cdot |$ denote the norm in $\X$ and
$\Y$ respectively. We have
\begin{equation*}
T(g, \epsilon)=T(G,\epsilon)\quad \forall\epsilon > 0.
\end{equation*}
This implies that $T(G, \epsilon)$ is relatively dense in $\R$ for
all $\epsilon > 0$, i.e. $G \in AP(\Y)$.

Next we prove that $\sigma_{b}(G)=\sigma_{b}(g)$. In fact, we have
\begin{eqnarray*}
a(\lambda, G)(s) &=& \lim_{T\rightarrow\infty}\frac{1}{2T}\int^{T}_{-T}e^{-i\lambda t}S(t)g(s)dt\\
 &=& \lim_{T\rightarrow\infty}\frac{1}{2T}\int^{T}_{-T}e^{-i\lambda t}g(t+s)dt\\
 &=& \lim_{T\rightarrow\infty}\frac{1}{2T}\int^{T+s}_{-T+s}e^{-i\lambda (\tau-s)}g(\tau)d\tau\\
 &=& e^{i\lambda s}a(\lambda, g),\quad \forall s \in \R .
\end{eqnarray*}
So $a(\lambda, G)=0$ if and only if $a(\lambda, g)=0$, i.e.
$\sigma_{b}(G)=\sigma_{b}(g)$.
\end{proof}

\begin{definition}\label{def 2-2}
A function $u\in BUC(\R,\X)$ is said to be a {\it strong mild
solution} of Eq.(\ref{1}) if there exists a sequence of classical
solutions $u_n$ to Eq.(\ref{1}) with $f$ replaced by $f_n\in
BUC(\R,\X)$ such that $\lim_{n\to \infty} u_n=u$, and $\lim_{n\to
\infty} f_n=f$ in the sup-norm topology of $BUC(\R,\X)$.
\end{definition}
We remark that this concept of strong mild solutions is given in
\cite{pru} as mild solutions. By the proof of Lemma \ref{lem 2} it
is clear that each strong mild solution of Eq.(\ref{1}) is a mild
solution of Eq.(\ref{1}). This justifies the terminology we use
here. In \cite{paz} the reader can find the equivalence of these two
concepts of mild solutions when $A$ generates a $C_0$-semigroup and
${\cal B}=0$.

\begin{lemma}\label{lem 6}
Let $u\in AP(\X)$ be a strong mild solution of Eq.(\ref{1}), $U(t)
:= S(t)u, F(t) := S(t)f. $ Then
\begin{equation}\label{3}
 a(\lambda, U)\in D(\mathcal{A}),  \quad \forall \lambda \in \R,
\end{equation}
\begin{equation}\label{4}
a(\lambda, u)\in D(A), \quad \forall \lambda \in \R,
\end{equation}
and
\begin{equation}\label{5}
(i\lambda - \mathcal{A} - \mathcal{B}) a(\lambda, U) = a(\lambda,
F).
\end{equation}
\end{lemma}
\begin{proof}
Since $u$ is a strong mild solution of Eq.(\ref{1}), there exist
$\{f_n\}\subset BUC(\R, \X)$ and a sequence of classical solutions
$\{u_n\}\subset BUC(\R, \X) $ such that $u_n \rightarrow u, f_n
\rightarrow f,$ and
$$
 \dot u_n(t) = Au_n(t) + [{\cal B}u_n](t)+ f_n(t).
$$
Letting $U_n(t) := S(t)u_n, F_n(t) := S(t)f_n, $ we define the
following
\begin{equation*}
a(\lambda, G; T)(s) = \frac{1}{2T}\int^{T}_{-T}e^{-i\lambda
t}S(t)g(s)dt,\quad \forall \lambda\in \R , T>0, s\in \R
\end{equation*}
for all $g\in BUC(\R,\X)$, where $G(t)$ is defined as Lemma \ref{lem
5}. Then
\begin{eqnarray*}
a(\lambda, S(t)\dot{u}_n; T)(s) &=& \frac{1}{2T}\int^{T}_{-T}e^{-i\lambda t}S(t)\dot{u}_ndt\\
 &=& \frac{1}{2T}\int^{T}_{-T}e^{-i\lambda t}\dot{u}_n(t+s)dt\\
 &=& \frac{1}{2T}\int^{T+s}_{-T+s}e^{-i\lambda (\tau-s)}\dot{u}_n(\tau)d\tau\\
 &=& i\lambda \frac{1}{2T}\int^{T+s}_{-T+s}e^{-i\lambda (\tau-s)}u_n(\tau)d\tau+
 \frac{1}{2T}[u_n(T+s)e^{-i\lambda T} - u_n(-T+s)e^{i\lambda T}]\\
 &=& i\lambda a(\lambda, U_n; T)(s)+
 \frac{1}{2T}[u_n(T+s)e^{-i\lambda T} - u_n(-T+s)e^{i\lambda T}].
\end{eqnarray*}
So we have
\begin{equation*}
(i\lambda - \mathcal{A} - \mathcal{B}) a(\lambda, U_n; T)(s) +
\frac{1}{2T}[u_n(T+s)e^{-i\lambda T} - u_n(-T+s)e^{i\lambda T}] =
a(\lambda, F_n; T)(s).
\end{equation*}
Now letting $n\rightarrow \infty$, by the closedness of $(i\lambda
- \mathcal{A} - \mathcal{B})$ we obtain
\begin{equation*}
(i\lambda - \mathcal{A} - \mathcal{B}) a(\lambda, U; T)(s) +
\frac{1}{2T}[u(T+s)e^{-i\lambda T} - u(-T+s)e^{i\lambda T}] =
a(\lambda, F; T)(s).
\end{equation*}
By Lemma \ref{lem 5}, $a(\lambda, U)$ and $a(\lambda, F)$ are well
defined, and
\begin{equation*}
\lim_{T\rightarrow\infty}a(\lambda, U; T)=a(\lambda, U),\quad
\lim_{T\rightarrow\infty}a(\lambda, F; T)=a(\lambda, F).
\end{equation*}
As $T\rightarrow \infty$, by the closedness of $(i\lambda -
\mathcal{A} - \mathcal{B})$ again we obtain (\ref{3}) and (\ref{5}).
In the proof of Lemma \ref{lem 5} we proved that $a(\lambda, U)(s) =
e^{i\lambda s}a(\lambda, u), \forall s \in \R $. Letting $s=0$, we
obtain (\ref{4}).
\end{proof}

The following theorem is an immediate consequence of Lemmas
\ref{lem 5} and \ref{lem 6}.
\begin{theorem}\label{th 6}
Let $u\in AP(\X)$ be a strong mild solution of Eq.(\ref{1}). Then
\begin{equation}\label{boh spe est}
\sigma_b(f) \subset \sigma_b(u) \subset \sigma_i(\mathcal{A +
B})\cup \sigma_b(f) .
\end{equation}
\end{theorem}
As consequences of Theorem \ref{th 6} we have.
\begin{corollary}\label{cor 10}
Let $u\in AP(\X)$ be a strong mild solution of Eq.(\ref{1}). Then
for any finite subset $\Lambda^1=\{\lambda_1, \lambda_2, \cdots,
\lambda_N\}\subset\sigma _i ({\cal A}+{\cal
B})\backslash\sigma_b(f)$ , there exists a mild solution $u_f\in
AP(\X) $, such that
\begin{equation}\label{9}
\sigma_b(u_f) \subset (\sigma_i(\mathcal{A + B})\backslash
\Lambda^1)\cup \sigma_b(f).
\end{equation}
\end{corollary}
\begin{proof}
Let $\Lambda:=\sigma_i(\mathcal{A + B})\cup \sigma_b(f),
\Lambda^2:=\Lambda\backslash \Lambda^1 = (\sigma_i(\mathcal{A +
B})\backslash \Lambda^1)\cup \sigma_b(f)$, and define
\begin{eqnarray*}
\Lambda_{AP} &:=& \{g\in AP(\X) : \sigma_b(g)\subset \Lambda\},\\
\Lambda^1_{AP} &:=& \{g\in AP(\X) : \sigma_b(g)\subset \Lambda^1\},\\
\Lambda^2_{AP} &:=& \{g\in AP(\X) : \sigma_b(g)\subset \Lambda^2\}.
\end{eqnarray*}
By the definition of Bohr spectrum,
$\Lambda_{AP},\Lambda^1_{AP},\Lambda^2_{AP}$ are closed linear
subspaces of $AP(\X)$ and $\Lambda^1_{AP}\bigcap\Lambda^2_{AP} =
\{0\}$. So we can easily obtain the decomposition
\begin{equation*}
\Lambda_{AP}=\Lambda^1_{AP}\oplus\Lambda^2_{AP}.
\end{equation*}
Therefore, there exist $u_1 \in \Lambda^1_{AP}, u_2\in
\Lambda^2_{AP}$ , where
\begin{equation*}
u_1 = \sum_{j=1}^{N}a(\lambda_j, u)e^{i \lambda_j t},
\end{equation*}
such that $u = u_1 + u_2$. By Lemma \ref{lem 6}, $a(\lambda_j, u)
\in D(A), j=1, 2, \cdots , N$. Letting $f_1 := (L_\Lambda -{\cal B})
u_1, f_2 := (L_\Lambda -{\cal B}) u_2$, we have
\begin{equation}\label{10}
(L_\Lambda -{\cal B}) u= f_1 + f_2 = f.
\end{equation}
Since $u_1$ is differentiable and $u_1(t)\in D(A)$ for all $t$ one
may check that $u_1$ is a classical solution of Eq. (\ref{1}) with
$f$ replaced by $f_1$. This yields that $u_2$ is a strong mild
solutions of Eq. (\ref{1}) with $f$ replaced by $f_2$. Now by
Theorem \ref{th 6}, we obtain $\sigma_b(f_1) \subset
\sigma_b(u_1)\subset\Lambda^1,
\sigma_b(f_2)\subset\sigma_b(u_2)\subset \Lambda^2$. This implies
that $f_1 \in \Lambda^1_{AP}, f_2 \in \Lambda^2_{AP}$. On the other
hand, $f\in \Lambda^2_{AP}$, by (\ref{10}) we have $f_1=0$ and
$f_2=f$. Therefore,
\begin{equation}\label{101}
(L_\Lambda -{\cal B}) u_2 = f.
\end{equation}
This shows that $u_f :=u_2$ is a mild solution of Eq.(\ref{1}) and
(\ref{9}) holds.
\end{proof}

\begin{corollary}\label{cor 11}
Let $u\in AP(\X)$ be a strong mild solution of Eq.(\ref{1}), and let
 $\sigma _i ({\cal A}+{\cal B})$ be finite. Then there exists a
mild solution $u_f\in AP(\X) $, such that $\sigma_b(u_f) =
\sigma_b(f)$.
\end{corollary}
\begin{proof}
It is obtained by Theorem \ref{th 6} and Corollary \ref{cor 10}
immediately.
\end{proof}

Consider the case where $f$ is a quasi-periodic function. We have
the following result on the non-existence of $k$-quasi-periodic
solutions.

\begin{corollary}\label{cor 12}
Let $f$ be a $l$-quasi periodic function. Then there exists no
k-quasi periodic strong mild solution of Eq.(\ref{1}) for any $k<l$.
\end{corollary}
\begin{proof}
In fact, if $u$ is $k$-quasi periodic strong mild solution, then by
 Theorem \ref{th 6} we have
$ \sigma_b (f) \subset \sigma_b(u) $. This implies that $f$ is
$k$-quasi periodic. Therefore, this is a contradiction.
\end{proof}

Another application of Theorem \ref{th 6} is illustrated in the
following
\begin{example}
 If $f\in AP (\X)$ is not $2\pi$-periodic, that is $\sigma_b(f)\not\subset 2\pi\Z$, then Eq.(\ref{1}) has no
 $2\pi$-periodic strong mild solutions.
\end{example}

\section{Applications and Examples}
Since our operator ${\cal B}$ in Eq.(\ref{1}) is very general, our
results in the previous section cover many classes of functional
differential equations that are considered in the literature. To
illustrate this, let $B(\cdot )$ be a function with bounded
variation that takes values in $L(\X)$, where $L(\X)$ denotes the
space of all bounded linear operators from $\X$ to $\X$, and let
\begin{equation}\label{b}
\int^\infty_{-\infty} \| dB(\eta ) \| e^{\delta | \eta|} < \infty
\end{equation}
for some $\delta >0$. Then we can define
$$
[{\cal B}u](t)= \int^\infty_{-\infty}dB(\eta ) u(t+\eta ), \quad
\forall u\in BUC(\R,\X), \ t\in \R .
$$
Consider equations of the form
\begin{equation}\label{22}
\dot u(t) = Au(t) +\int^\infty _{-\infty} dB(\eta ) u(t+\eta ) +
f(t),\quad t\in \R .
\end{equation}
Let $u\in BUC(\R,\X)$ be a mild solution of (\ref{22}). By the same
method as in the proof of Lemma \ref{lem 3} we can show that ${\cal
D}-{\cal A}-{\cal B}$ is closable in $BUC(\R,\X)$. Moreover, because
of the closedness of the operator $L$ with domain consisting of all
$u\in BUC(\R,\X)$ such that $\int^t_s u(\xi )d\xi \in D(A)$ for all
$t\ge s$, and
$$
u(t)-u(s) = A\int^t_s u(\xi )d\xi +\int^t_s g(\xi )d\xi ,
$$
for some $g\in BUC(\R,\X)$, and $Lu=g$, the operator $L-{\cal B}$ is
a closed extension of ${\cal D}-{\cal A}-{\cal B}$. This implies
that each mild solution in the sense of Definition \ref{def 2}
satisfies the equation
\begin{equation}
u(t)-u(s) = A\int^t_s u(\xi )d\xi + \int^t_s \int^\infty_{-\infty}
dB(\eta )u(\eta +\xi )d\xi +\int^t_sf(\xi )d\xi ,\quad t\ge s.
\end{equation}
In this case we will show that the set $\sigma_i(\mathcal{A+B})$
can be replaced by a simpler one that can be easier computed. In
the proof of Theorem \ref{th 1} we got the following equation
\begin{equation*}
(\lambda-\mathcal{D})^{-1}f+u=[\lambda-(\mathcal{A+B})](\lambda-\mathcal{D})^{-1}u,\quad
Re\lambda\neq0.
\end{equation*}
Recall that
\begin{eqnarray*}
(\lambda -{\cal D})^{-1}u (\xi )=
\begin{cases}
\begin{array}{ll} \int^  {\infty}_0  e^{-\lambda t }u(t+\xi
)d\xi,
&(\mbox{if}\ Re\lambda > 0),\\ \\
-\int _ {-\infty}^0  e^{-\lambda t}f(t+\xi )d\xi,     & (\mbox{if }
\ Re\lambda < 0)
\end{array} \end{cases}
\end{eqnarray*}
and
\begin{equation*}
\hat{u} (\lambda ) = (\lambda -{\cal D})^{-1}u (0 ).
\end{equation*}
We have
\begin{eqnarray*}
(\lambda-\mathcal{D})^{-1}f(0)+u(0)&=&[\lambda-(\mathcal{A+B})](\lambda-\mathcal{D})^{-1}u(0)\\
 &=& [\lambda(\lambda-\mathcal{D})^{-1}u-\mathcal{A}(\lambda-\mathcal{D})^{-1}u-(\lambda-\mathcal{D})^{-1}\mathcal{B}u](0)\\
 &=& \lambda
 \hat{u}(\lambda)-A\hat{u}(\lambda)-[(\lambda-\mathcal{D})^{-1}\mathcal{B}u](0).
\end{eqnarray*}
In the case of $\delta > Re \lambda>0$ we have
\begin{eqnarray*}
[(\lambda-\mathcal{D})^{-1}\mathcal{B}u](0)&=&\int^  {\infty}_0
e^{-\lambda t}[\mathcal{B}u](t)dt\\
&=& \int^  {\infty}_0
e^{-\lambda t}\int^\infty_{-\infty} dB(\eta )u(\eta +t )dt\\
&=& \int^\infty_{-\infty} dB(\eta )\int^{\infty}_0 e^{-\lambda
t}u(\eta +t )dt\\
&=& \int^\infty_{-\infty} dB(\eta )e^{-\lambda \eta}(\int^{\infty}_0
e^{-\lambda \xi}u(\xi)d\xi - \int^{\eta}_0 e^{-\lambda
\xi}u(\xi)d\xi )\\
&=&\int^\infty_{-\infty} dB(\eta )e^{-\lambda \eta}\hat{u}(\lambda)
- \int^\infty_{-\infty} dB(\eta )e^{-\lambda \eta}\int^{\eta}_0
e^{-\lambda \xi}u(\xi)d\xi.
\end{eqnarray*}
In the case of $-\delta < Re \lambda<0$, we obtain a similar
equality. So we have
\begin{equation*}
\hat{f} (\lambda ) + u(0) = \lambda \hat{u}(\lambda)
-A\hat{u}(\lambda)-\int^\infty_{-\infty} dB(\eta )e^{-\lambda
\eta}\hat{u}(\lambda)+\int^\infty_{-\infty} dB(\eta )e^{-\lambda
\eta}\int^{\eta}_0 e^{-\lambda \xi}u(\xi)d\xi ;
\end{equation*}
that is,
\begin{equation*}
\hat{f} (\lambda ) + u(0)-\int^\infty_{-\infty} dB(\eta )e^{-\lambda
\eta}\int^{\eta}_0 e^{-\lambda \xi}u(\xi)d\xi= (\lambda
 -A-\int^\infty_{-\infty} dB(\eta
)e^{-\lambda \eta})\hat{u}(\lambda).
\end{equation*}
If $\xi\not\in \Sigma$, where
\begin{equation}
\Sigma :=\{z\in\C: \ |Re z|<\delta ,  \not\exists \ \left(z
 -A-\int^\infty_{-\infty} dB(\eta)e^{-z \eta}\right)^{-1}\in L(\X) \} ,
\end{equation}
we have
\begin{equation*}
\hat{u}(\lambda)=(\lambda -A-\int^\infty_{-\infty} dB(\eta
)e^{-\lambda \eta})^{-1}(\hat{f} (\lambda ) +
u(0)-\int^\infty_{-\infty} dB(\eta )e^{-\lambda \eta}\int^{\eta}_0
e^{-\lambda \xi}u(\xi)d\xi) .
\end{equation*}
We can prove (see e.g. \cite{furnaimin}) that if $i\xi\not\in
\Sigma$ and $\xi\not\in sp(f)$, then
$$
\left(\lambda -A-\int^\infty_{-\infty} dB(\eta )e^{-\lambda
\eta}\right)^{-1}, \quad \int^\infty_{-\infty} dB(\eta
)e^{-\lambda \eta}\int^{\eta}_0 e^{-\lambda \xi}u(\xi)d\xi
$$
are analytic in $\lambda$ around $i\xi$. This implies $\xi\not\in
sp(u)$, so $sp(u)\subset\Sigma_i\cup sp(f)$, where $\Sigma_i:=\{
\xi \in \Sigma :\ i\xi \in \Sigma\} $. So we arrive at
\begin{proposition}
Let the above mentioned conditions for Eq. (\ref{22}) be
satisfied, and let $u\in BUC(\R,\X)$ be a mild solution of Eq.
(\ref{22}). Then,
\begin{equation}\label{4.4}
sp( u) \subset \Sigma_i \cup sp(f) .
\end{equation}
\end{proposition}
The advantage of the estimate (\ref{4.4}) is that we need only to
study the "{\it characteristic roots}" of the equation, instead of
the spectrum $\sigma_i ({\cal A}+{\cal B})$. We are ready to see
how our results obtained above extend respective ones in
\cite{furnaimin}. Now the reader can re-state all above conditions
and results in terms of the spectral set $\Sigma_i$, instead of
$\sigma_i ({\cal A}+{\cal B})$.

\end{document}